\documentclass{article}
\usepackage[dvips]{graphicx}
\usepackage{amssymb}
\usepackage{amsmath}

\usepackage{caption}
\usepackage{hyperref}
\usepackage{arcs}
\usepackage {xcolor} 
\usepackage{subfigure}
\usepackage{authblk}
\usepackage{calrsfs}

\newtheorem{theorem}{Theorem}

\newtheorem{Corollary}{Corollary}

\usepackage{graphicx} % Required for inserting images

\newenvironment{AMS}{\small\bf 2020 AMS subject classification: }{}

\title{Unisolvence of unsymmetric random Kansa collocation 
by Gaussians and other analytic RBF vanishing at infinity}
\author{Alvise Sommariva and Marco Vianello\\University of Padova}
\date{May 2024}

\begin{document}

\maketitle

\begin{abstract}
We give a short proof of almost sure invertibility of unsymmetric random Kansa collocation matrices by a class of analytic RBF vanishing at infinity, for the Poisson equation with Dirichlet boundary conditions. Such a class includes 
popular Positive Definite instances such as Gaussians, Generalized Inverse MultiQuadrics and Mat\'ern RBF. The proof works on general domains in any dimension, with any distribution of boundary collocation points and any continuous random distribution of internal collocation points.
\end{abstract}

\vskip0.5cm
\noindent
\begin{AMS}
{\rm 65D12, 65N35.}
\end{AMS}
\vskip1cm

\section{Introduction} 
The unsymmetric Kansa method for meshless discretization of PDEs, in both the square (standard collocation) and rectangular (overtesting) formulations, has a long and successful history in applications; see, with no pretence of exhaustivity,  \cite{CL20,CFC14,CLS18,F07,K86,LOS06,SW06,W05} with the references therein. 
Nevertheless, existence of sufficient conditions ensuring invertibility of unsymmetric Kansa collocation matrices is still a substantially open problem, as it has been for more than 20 years after the breakthrough paper by Hon and Schaback \cite{HS01}, who showed that there exist point configurations leading to singularity of the collocation matrices (though these are very special and ``rare''cases).
The lack of clear well-posedness conditions for the discrete problem has been considered one of the main drawbacks of unsymmetric Kansa collocation, despite its manifest effectiveness in many applications.

In two recent papers, a new approach to prove unisolvence in the case of {\em random} collocation has been proposed, for Thin-Plate Splines, MultiQuadrics and Inverse MultiQuadrics. The key property is that each basis function is analytic but presents some singularity either in the real (TPS,\cite{DASV24}) or in the complex domain (MQ and  IMQ, \cite{CDRDASV24}), which is not a singularity for the other basis elements. A similar approach has been used also for polynomial-free random interpolation by general Polyharmonic Splines \cite{BSV24,SV24} and Generalized Multiquadrics \cite{SV24-2}. However, it does not work for example in the case of Gaussians, which are entire functions. 

In this paper we study unisolvence of unsymmetric random Kansa collocation for the Poisson equation by a class of analytic RBF with scaled radial function
$$
\phi_\varepsilon(r)=\phi(\varepsilon r)
$$
which is $C^2$ in $[0,+\infty)$ and analytic in $(0,+\infty)$, such that each RBF and its Laplacian vanish at $\infty$ but the Laplacian does not vanish at the RBF center. 
As we shall see, such a class includes 
popular Positive Definite instances such as Gaussians, Generalized Inverse MultiQuadrics and Mat\'ern RBF. 
The scale $\varepsilon>0$ represents the so-called \emph{shape parameter} associated with RBF and as known is relevant in order to control
the trade-off between conditioning and accuracy; cf. e.g. \cite{F07,LS23}. 

We consider in particular the Poisson equation with Dirichlet boundary conditions (cf. e.g. \cite{E98})
\begin{equation}\label{poisson}
\left\{
\begin{array}{l}
\Delta u(P)=f(P)\;,\;P\in \Omega\;,
\\
u(P)=g(P)\;,\;P\in \partial \Omega\;,
\end{array} 
\right .
\end{equation}
where $\Omega\subset \mathbb{R}^d$ is a bounded domain (connected open set), $P=(x_1,\dots,x_d)$ and $\Delta=\partial ^2/\partial x_1^2+\dots+\partial ^2/\partial x_d^2$ is the Laplacian. In \cite{DASV24}, where we studied Kansa discretization by Thin-Plate Splines in $\mathbb{R}^2$, we assumed that the boundary is a curve possessing an analytic parametrization. Differently here, as in \cite{CDRDASV24} (which concerns MultiQuadrics and Inverse MultiQuadrics), we do not make any restrictive assumption on $\partial\Omega$, except for the usual ones that guarantee well-posedness and regularity of the solution (cf. e.g. \cite{E98}). 

We recall that unsymmetric Kansa collocation  consists in seeking a function
\begin{equation} \label{u_N}
u_N(P)=\sum_{j=1}^n{c_j\,\phi_j(P)}+\sum_{k=1}^m{d_k\,\psi_k(P)}\;,\;\;N=n+m\;,
\end{equation}
where
\begin{equation} \label{phij}
\phi_j(P)=\phi_\varepsilon(\|P-P_j\|_2)\;,\;\;\{P_1,\dots,P_n\}\subset \Omega\;,
\end{equation}
\begin{equation} \label{psik}
\psi_k(P)=\phi_\varepsilon(\|P-Q_k\|_2)\;,\;\;\{Q_1,\dots,Q_m\}\subset 
\partial\Omega\;,
\end{equation}
such that
\begin{equation}\label{poisson_disc}
\left\{
\begin{array}{l}
\Delta u_N(P_i)=f(P_i)\;,\;i=1,\ldots,n
\\
u_N(Q_h)=g(Q_h)\;,\;h=1,\ldots,m\;.
\end{array} 
\right .
\end{equation}

Kansa collocation can be rewritten in matrix form as
\begin{equation} \label{Kansa-system}
\left(\begin{array}{cc}
\Delta\Phi & \Delta\Psi\\
\\
\Phi & \Psi
\end{array} \right) 
\left(\begin{array}{c}
\mathbf{c}\\
\\
\mathbf{d}
\end{array} \right) 
=
\left(\begin{array}{c}
\mathbf{f}\\
\\
\mathbf{g}
\end{array} \right) 
\end{equation}
where the $N\times N$ block matrix is 
$$
K_N=K_N(\{P_i\},\{Q_h\})=\left(\begin{array}{cc}
\Delta\Phi & \Delta\Psi\\
\\
\Phi & \Psi
\end{array} \right) 
=\left(\begin{array}{cc}
(\Delta\phi_j(P_i)) & (\Delta\psi_k(P_i))\\
\\
(\phi_j(Q_h)) & (\psi_k(Q_h))
\end{array} \right) 
$$
and $\textbf{f}=\{f(P_i)\}$, $\textbf{g}=\{g(Q_h)\}$, $1\leq i,j\leq n$, $1\leq h,k\leq m$.

\section{Unisolvence of random Kansa collocation}

We can now state and prove our main result on unisolvence of random Kansa collocation by a class of analytic RBF vanishing at infinity. Thus will be followed by a Corollary on the application to Gaussians and other Positive Definite RBF. 
In the sequel we shall denote by $H(D)$ the space of analytic functions in the open connected set $D\subset \mathbb{R}^d$. 

\begin{theorem}
Let $\phi:[0,+\infty)\to \mathbb{R}$ be a radial function such that:
\begin{itemize}
\item[$(i)$] $\phi\in C^2([0,+\infty))\cap H((0,+\infty))$, $\lim_{r\to \infty}\phi(r)=0$;
\item[$(ii)$] $\ell(r)=\phi''(r)+\phi'(r)/r$ is continuous at $r=0$, 
$\ell(0)\neq 0$, $\lim_{r\to \infty}\ell(r)=0$;
\item[$(iii)$] the RBF interpolation matrix $V_m=\phi_\varepsilon(\|Q_h-Q_k\|)$, $1\leq h,k\leq m$, is nonsingular for every set of distinct points 
$\{Q_1,\dots,Q_m\}\subset \mathbb{R}^d$.
\end{itemize}
Moreover, $K_n$ be the Kansa collocation matrix by $\phi_\varepsilon(r)=\phi(\varepsilon r)$ defined above for equation (\ref{poisson}) in $\Omega\subset \mathbb{R}^ d$, $d\geq 2$, where $\{Q_h\}$ is any fixed set of $m$ distinct points on $\partial\Omega$, and  $\{P_i\}$ is a sequence of i.i.d. (independent and identically distributed) random points in $\Omega$ with respect to any probability density $\sigma \in L^1_+(\Omega)$. 

Then for every $m \geq 1$ and for every $n\geq 0$ the matrix $K_n$ is a.s. (almost surely) nonsingular.
\end{theorem}
\vskip0.2cm 
\noindent 
{\bf Proof.} Let us fix $\varepsilon>0$ and define $\phi_A(P)=\phi_\varepsilon(\|P-A\|_2)$.
Then we have $\phi_A(B)=\phi_B(A)$ and $\phi_A(A)=\phi(0)$. Moreover (cf. \cite[Appendix D]{F07})
$$
\Delta\phi_A(P)=\varepsilon^ 2(\phi''(\varepsilon\|P-A\|_2)+\phi'(\varepsilon\|P-A\|_2)/(\varepsilon \|P-A\|_2))=\varepsilon^2\ell(\varepsilon \|P-A\|_2)
$$
so that $\Delta\phi_A(B)=\Delta\phi_B(A)$ and $\Delta\phi_A(A)=\varepsilon^2\ell(0)$. Observe in addition that $\phi_A(\cdot), \Delta\phi_A(\cdot)\in H(\mathbb{R}^d\setminus \{A\})$ by analiticity of $\phi$ and $\sqrt{\cdot}$ in $\mathbb{R}^+$, and by $(i),(ii)$ 
\begin{equation} \label{limit}
\lim_{P\to \infty}\phi_A(P)=0=\lim_{P\to \infty}\Delta \phi_A(P)\;.
\end{equation}

The proof proceeds by induction on $n$. The assertion is certainly true for $n=0$ by $(iii)$, since $K_0=V_m$. Let us assume that the assertion is true for a fixed $n\geq 0$ and define the matrix
$$
K(P)=
\left(\begin{array}{ccc}
\Delta\Phi  & \Delta\Psi & (\Delta\vec{\phi}(P))^t\\
\\

\Phi & \Psi & (\vec{\psi}(P))^t \\ \\
\Delta\vec{\phi}(P) &\Delta\vec{\psi}(P) & \varepsilon^2 \ell(0)
\end{array} \right) 
$$
where $\vec{\phi}(P)=(\phi_1(P),\dots,\phi_{n}(P))$ and $\vec{\psi}(P)=(\psi_1(P),\dots,\psi_{m}(P))$.
Observe that $|\mbox{det}(K_{n+1})|=|\mbox{det}(K(P_{n+1}))|$, since  $\psi_k(P_{n+1})=\phi_{n+1}(Q_k)$,  
$\Delta\phi_j(P_i)=\Delta\phi_i(P_j)$, and thus $K(P_{n+1})$ corresponds to a permutation of rows and columns of $K_{n+1}$. Now in view of (\ref{limit}) we have that
$$
\lim_{P\to \infty}K(P)=L
=\left(\begin{array}{cc}
K_n & \vec{0}^{\,t}\\
\\
\vec{0} & \varepsilon^2\ell(0)
\end{array} \right) 
$$
where the limit matrix $L$ is a.s. invertible, since $\mbox{det}(L)=\varepsilon^2\ell(0)\mbox{det}(K_n)$ and a.s. $\mbox{det}(K_n)\neq 0$ by inductive hypothesis while $\ell(0)\neq 0$ by $(ii)$.

This means that the function 
$$
F(P)=\mbox{det}(K(P))
$$
which is in $H(\mathbb{R}^d\setminus \{P_1,\dots,P_n\})$, a.s. does not vanish in a neighborhood of $\infty$ in $\mathbb{R}^d$, and hence is not identically zero in the open connected set $\mathbb{R}^d\setminus \{P_1,\dots,P_n\}$. 

In view of basic result in the theory of analytic functions,  
the zero set of $F$ in the open connected set $\mathbb{R}^d\setminus \{P_1,\dots,P_n\}$, say $Z(F)$, is then a null set, i.e. it has zero Lebesgue measure (cf. \cite{M20} for an elementary proof). Consequently, the zero set of $F$ in $\Omega$, say $Z_\Omega(F)$ is also a null set since 
$$
Z_\Omega(F)=(Z_\Omega(F)\cap \{P_1,\dots,P_n\})\cup 
(Z(F)\cap \Omega)
$$
where the first is a finite set and the second a subset of a null set, and hence is a null set also for $\sigma(x)dx$, i.e. the probability that $P_{n+1}$ falls in $Z_\Omega(F)$ is zero. 

To conclude, considering the probability of the corresponding events
$$ 
\mbox{prob}\{\mbox{det}(K_{n+1})=0\}=\mbox{prob}\{F(P_{n+1})=0\}
$$
$$
=\mbox{prob}\{F\equiv 0\}
+\mbox{prob}\{F\not\equiv 0\;\&\;P_{n+1}\in Z_\Omega(F)\}
=0+0=0\;,
$$ 
and the inductive step is completed.
\hspace{0.2cm} $\square$

\vskip0.5cm

\begin{Corollary}
Under the assumptions of Theorem 1 on the boundary and the internal collocation points, the Kansa collocation matrix for equation (\ref{poisson}) by Gaussian, Generalized Inverse MultiQuadric and Mat\'ern RBF is almost surely nonsingular.  
\end{Corollary}
\vskip0.2cm     
\noindent 
{\bf Proof.} It is sufficient to check properties $(i)-(iii)$ of Theorem 1. For all instances $(iii)$ holds since $\phi(\|x\|_2)$ is a strictly Positive Definite RBF and $V_m$ is a positive definite matrix (cf. e.g. \cite{F07}).

\begin{itemize}
\item Gaussians: $\phi(r)=e^{-r^2}$ and $\ell(r)=4e^{-r^2}(r^2-1)$ are entire functions which vanish at $+\infty$, $\ell(0)=4$;

\item Generalized Inverse MultiQuadrics: $\phi(r)=(1+r^2)^\beta$, $\beta<0$, and $\ell(r)=4\beta(1+r^2)^{\beta-2}(1+\beta r^2)$ are analytic in $\mathbb{R}$ and vanish at $\pm \infty$, $\ell(0)=4\beta$; this class includes the standard IMQ ($\beta=-1/2$) and the IQ (Inverse Quadratic, $\beta=-1$);

\item Mat\'ern: $\phi(r)=\frac{2^{1-\nu}}{\Gamma(\nu)}\,r^\nu K_\nu(r)$, where $K_\nu$ is the modified Bessel function of the second kind or Macdonald function of order $\nu>0$ (cf. e.g. \cite{NIST}). These  radial functions has been used for Kansa collocation e.g. in \cite{MB02} with half-integer parameter $\nu=k+1/2$ with $k\in \mathbb{N}$. Now, let us take in general $\nu \in \mathbb{R}$, $\nu>1$. It is known that $K_\nu$ is analytic in $\mathbb{C}\setminus (-\infty,0]$, and $K_\nu(r)\sim \sqrt{\pi/(2r)}\,e^{-r}$ as $r\to +\infty$ 
while $r^\nu K_\nu(r)\sim 2^{\nu-1}\Gamma(\nu)$ as $r\to 0^+$.
Moreover $(r^\nu K_\nu(r))'=-r^{\nu} K_{\nu-1}(r)$ so that $\phi'(r)=-\frac{2^{1-\nu}}{\Gamma(\nu)}\,r^\nu K_{\nu-1}(r)$ and 
$\phi''(r)=-\frac{2^{1-\nu}}{\Gamma(\nu)}(r \cdot r^{\nu-1}K_{\nu-1}(r))^{\prime}=-\frac{2^{1-\nu}}{\Gamma(\nu)}\,(r^{\nu-1}K_{\nu-1}(r)-r^{\nu}K_{\nu-2}(r))$, from which one easily gets  
$$
\ell(r)=-\frac{2^{1-\nu}}{\Gamma(\nu)}\,
(2r^{\nu-1}K_{\nu-1}(r)-r^{\nu}K_{\nu-2}(r))\;.
$$
From the relations above we can finally compute  
$\lim_{r\to +\infty}\phi(r)=0=\lim_{r\to +\infty}\ell(r)$ and $\lim_{r\to 0^+}\ell(r)=-\frac{\Gamma(\nu-1)}{\Gamma(\nu)}=\frac{1}{1-\nu}\neq 0\;$.
\hspace{0.2cm} $\square$

\end{itemize}

\vskip0.2cm

To our knowledge, this is the first set of sufficient conditions for unisolvence of unsymmetric Kansa collocation by Gaussian and Mat\'ern RBF. In the case of classical Inverse MultiQuadrics ($\beta=-1/2$) an alternative but more complicated unisolvence proof for random collocation has been recently provided in \cite{CDRDASV24}.

\vskip0.5cm 
\noindent
{\bf Acknowledgements.} 
\vskip0.25cm 

Work partially supported by the DOR funds of the University of Padova, and by the INdAM-GNCS 2024 Projects “Kernel and polynomial methods for approximation and integration: theory and application software''. 

This research has been accomplished within the RITA ``Research ITalian network on Approximation" and the SIMAI Activity Group ANA\&A, and the UMI Group TAA ``Approximation Theory and Applications".

\end{document}